\documentclass[12pt]{article}
\usepackage[margin=1in]{geometry}

\usepackage{url,subfig,mathrsfs,citesort,paralist}
\RequirePackage{amsthm,amsmath,amssymb}

\usepackage{color}

\theoremstyle{plain}
\newtheorem{theorem}{Theorem}

\theoremstyle{definition}

\renewcommand{\le}{\leqslant}
\renewcommand{\leq}{\leqslant}
\renewcommand{\ge}{\geqslant}
\renewcommand{\geq}{\geqslant}
\renewcommand{\P}{{\rm P}}
\renewcommand{\d}{{\rm d}}

\newcommand{\EE}{\mathrm{E}}
\newcommand{\R}{\mathbb{R}}
\newcommand{\RR}{\mathbb{R}}

\title{Delocalization of a $(1+1)$-dimensional\\stochastic wave 
	equation\footnote{Research supported in part by the NSF grants DMS-1307470.}
 }
\author{Jingyu Huang\footnote{Dept of Mathematics,
	University of Utah, Salt Lake City, UT 84112--0090. \texttt{jhuang@math.utah.edu}}
	\and
	Davar Khoshnevisan\footnote{Dept of Mathematics,
	University of Utah, Salt Lake City, UT 84112--0090. \texttt{davar@math.utah.edu}}}
\date{Version: October 24, 2016}

\begin{document}
\maketitle
\begin{abstract}
	A noteworthy property of many parabolic stochastic PDEs is that
	they locally linearize \cite{FKM,Hairer,HP,KSXZ}. We prove that, by contrast, 
	a large family of stochastic wave equations in dimension one do
	not possess this important property.\\


\end{abstract}

\section{Introduction}
Consider the following stochastic partial differential equation 
(SPDE, for short), indexed by space-time, $(0\,,\infty)\times\RR$:
\begin{equation}\label{eq: SHE}
	\partial_t v = \partial^2_{xx} v + \sigma(v)\xi
	\quad\text{subject to}\quad v(0) \equiv 1;
\end{equation}
where $v=v(t\,,x)$ for all space-time points
$(t\,,x)\in(0\,,\infty)\times\RR$; the forcing term
$\xi$ denotes space-time white noise; and $\sigma:\RR\to\RR$
is  a nonrandom, Lipschitz continuous function. It is well known in this, and related,
contexts, that the solution to \eqref{eq: SHE} locally linearizes. Indeed, let 
$Z=Z(t\,,x)$ denote the linearization of $v$; that is,
$Z$ solves the {SPDE} \eqref{eq: SHE} with $\sigma\equiv 1$. Then,
\[
	v(t\,, x+\varepsilon) - v(t\,,x) = \sigma(v(t\,,x)) 
	(Z(t\,,x+\varepsilon)-Z(t\,,x)) + \mathcal{R}_{t,x}(\varepsilon),
\]
where, as $\varepsilon\downarrow0$,  the remainder term $\mathcal{R}_{t,x}(\varepsilon)$
tends to zero much
faster than $Z(t\,,x+\varepsilon)-Z(t\,,x)$ does.
This fact appears explicitly---in different senses 
and contexts---in \cite{FKM,KSXZ,HP},
and in a different setting earlier in the fundamental works of Hairer \cite{Hairer,Hairer1}
respectively on the KPZ equation and on Hairer's theory of regularity structures.

The purpose of this short note is to point out
that, in sharp contrast with the parabolic setting, typical hyperbolic {SPDEs} 
do not locally linearize. To be concrete, let us consider a hyperbolic
SPDE of the type,
\begin{equation}\label{eq:SWE}
	\partial^2_{tt}u=  {\partial^2_{xx}u} + \sigma(u)\xi
	\quad\text{on $(0\,,\infty)\times\RR$,}
\end{equation}
subject to the initial conditions $u(0)=0$ and $(\partial_t u)(0)=1$,
to be concrete. Let $Y$ denote the linerization of $u$; that is, 
let $Y$ denote the solution to \eqref{eq:SWE} with $\sigma\equiv1$.

Let us suppose, to the contrary, that $u$ locally linearizes; that is,
let us posit that
\[
	u(t\,, x+h) - u(t\,,x) = \sigma(u(t\,,x)) \left(Y(t\,, x+h) - Y(t\,,x) \right) 
	+ \widetilde{\mathcal{R}}_{t,x}(h),
\]
where the remainder term {$\widetilde{\mathcal{R}}_{t,x}(h)$} is 
significantly smaller than $Y(t\,,x+h)-Y(t\,,x)$ when $h\approx 0$. Then,
a simple heuristic argument would suggest that, for every fixed $t>0$, the quadratic variation
of $[x_1\,,x_2]\ni x\mapsto u(t\,,x)$ would have to be a function of $u(t)$ alone;
in fact, a more careful heuristic
analysis of the random field $Y$ might suggest that the quadratic variation
of $u(t)$ is likely to be equal to
$2t \int_{x_1}^{x_2} [\sigma(u(t\,,x))]^2\,\d x$ at every fixed time point $t>0$.
Theorem \ref{thm:quadratic} refutes this assertion, and hence
rules out the possibility of the existence of good
local linearizations to $u$. 

\begin{theorem}\label{thm:quadratic}
	Choose and fix four real numbers $t>0$ and
	$x,X_1,X_2\in\R$ such that $X_1<x<X_2$.
	For all {integers} $N\ge1$ and  $1\le i\le N$, define
	$t_i =t_{i,N} := (i-1)t/N$ and $x_i = x_{i,N}:=
	X_1+ (i-1)(X_2-X_1)/N$ for integers $1 \leq i \leq N$. Then,
	the following is valid for every real number $p\ge2$.
	As $N\to\infty$:
	\begin{align}
		&\sum_{i=1}^N 
			\left[u(t_{i+1}\,, x) - u(t_i\,, x) \right]^2  
			\xrightarrow{\text{in $L^p(\Omega)$}} 
			\int_{Q(x,t)} \left[\sigma(u(s-|x-y|\,,y))\right]^2\,\d s\,\d y;
			\text{ and}
			\label{eq:time quadratic}\\
		&\sum_{i=1}^N 
			\left[ u(t\,,x_{i+1}) - u(t\,,x_i) \right]^2 
			\xrightarrow{\text{in $L^p(\Omega)$}}  
			\int_{(-t,t)\times(X_1,X_2)}
			\left[\sigma(u(s\,,y-|t-s|))\right]^2 \,\d s\,\d y;
			\label{eq:space quadratic}
	\end{align}
	where $Q(x\,,t):=\{(s\,,y): |y-x|\le|t-s|\}$.
\end{theorem}

The proof of \eqref{eq:time quadratic} hinges on the observation
that $u(r\,,y)\approx u(t_i - |x-y|\,, y)$
uniformly for all space-time points 
\begin{equation}\label{rQQ}
	(r\,,y)\in Q(x\,,t_{i+1}\,,t_i) := Q(x\,, t_{i+1})\setminus Q(x\,,t_{i});
\end{equation}
see Section \ref{section:Proof of thm quadratic}. 
This kind of approximation procedure also yields the following limit theorem
as a by product. In order to state the next result, 
first let $\mathcal{G}^0_{t,x}$ denote the $\P$-completion of the
$\sigma$-algebra generated by
all Wiener integrals of the form $\int\phi\,\d\xi$,
where $\phi$ is smooth and support in $Q(x\,,t)$, and $(t\,,x)\in\RR_+\times\RR$.
Then, define 
\[
	\mathcal{G}_{t,x} := \bigcap_{s > t} \mathcal{G}^0_{s,x}
	\qquad\text{for every $(t\,,x)\in\RR_+\times\RR$}.
\]

\begin{theorem}\label{thm:CLT and LIL}
	Choose and fix $t>0$ and $x\in\R$. Then,
	\begin{equation}\label{eq:CLT}
		 \frac{u(t+\varepsilon\,, x) - u(t\,,x) }{\sqrt\varepsilon} 
		 \xrightarrow{d}\int_{x-t}^{x+t} \sigma(u(t-|x-y|\,, y))\,W(\d y)
		 \qquad\text{as $\varepsilon\downarrow0$},
	\end{equation} 
	where $W$ is a standard two-sided Brownian motion that is
	independent of $\mathcal{G}_{t,x}$. Moreover,
	\begin{equation}\label{eq:LIL}
		\limsup_{\varepsilon\downarrow 0} 
		\frac{u(t+\varepsilon\,, x)-u(t\,,x)}{\sqrt{2\varepsilon\log \log(1/\varepsilon)}} = 
		\sqrt{\int_{x-t}^{x+t} \left[\sigma(u(t-|x-y|\,, y))\right]^2\,\d y} 
		\qquad\text{a.s.}
	\end{equation}
\end{theorem}

Theorems \ref{thm:quadratic} and \ref{thm:CLT and LIL}
are proved respectively in Sections \ref{section:Proof of thm quadratic} 
and \ref{sec:Proof of thm CLT and LIL}. 

\section{Proof of Theorem \ref{thm:quadratic}}\label{section:Proof of thm quadratic}

First let us recall that the random field 
$u = \{u(t\,,x)\}_{t\geq 0, x \in \RR}$ is a mild solution to \eqref{eq:SWE};
that is:
\begin{compactitem}
\item  $\{u(t)\}_{t\ge0}$ is predictable with respect to the Brownian filtration
	generated by  all Wiener integrals of the form $\int_{[0,t]\times\RR}\phi\,\d\xi$
	for  $\phi\in L^2(\RR_+\times\RR)$; and
\item For each $t\geq 0$ and $x \in \RR$, 
	\begin{equation}\label{mild}
		u(t\,, x )= 1 + \int_{Q(x,t)} \sigma(u(s\,,y))
		\,\xi(\d s\,\d y)\qquad\text{a.s.}, 
	\end{equation}
	where $Q(x\,,t)$ is defined in Theorem \ref{thm:quadratic},
	 and the stochastic integral is understood in the sense of Walsh \cite{Walsh}.
\end{compactitem}
It is well known that
the methods of Walsh's theory imply that \eqref{eq:SWE} is well posed
in the sense that there exists a unique continuous solution to the integral 
equation \eqref{eq:SWE}. 
It is also well known that for every $k\ge 2$,
\begin{equation}\label{thm:Holder continuity}
	\EE \left( |u(t_1\,, x_1)- u(t_2\,, x_2)|^k\right) \lesssim
	 |t_1-t_2|^{k/2} + |x_1-x_2|^{k/2},
\end{equation}
uniformly for all $t_1,t_2\in[0\,,t]$ and $x_1,x_2\in\R$;
consult Peszat and Zabczyk \cite{PZ}. This result appears, implicitly, earlier
in Dalang and Frangos \cite{DF}. Moreover, nearly all of the preceding,
and much more, is 
included in the  general theories of Dalang \cite{Dalang}
and Peszat and Zabzcyk \cite{PZ}.

We will prove \eqref{eq:time quadratic} and \eqref{eq:space quadratic} in  
successive steps, and in this order. \\

\emph{Step 1.} 
Recall the sets $Q(x\,,t_{i+1}\,,t_i)$, defined earlier in \eqref{rQQ}, and set
\[
	\mathcal{A}_N(x) := \sum_{i=1}^N \left[ u(t_{i+1}\,,x)-u(t_i\,,x) \right]^2\quad\text{\&}\quad
	\mathcal{B}_N(x) :=\sum_{i=1}^N \left(\int_{Q(x, t_{i+1}, t_i)} 
		\sigma \left(u(r_i(y)\,,y) \right)\xi(\d s\,\d y) \right)^2,
\]
where $r_i(y) :=\max( t_i-|x-y|\,, 0)$. 

In this first step of the proof, we will prove that $\mathcal{A}_N(x)-\mathcal{B}_N(x)\to0$
in $L^p(\Omega)$ as $N\to\infty$. 
To this end, let us first observe that,
for all $1\leq i \leq N$,  the mild formulation of the solution $u$
yields the  representation,
$u(t_{i+1}\,, x) - u(t_i\,, x) = \int_{Q(x, t_{i+1}, t_i)} 
\sigma(u(s\,,y)) \,\xi(\d s\,\d y).$
Therefore, if we write $\|\,\cdots\|_p$ for the $L^p(\Omega)$-norm, then
$\|\mathcal{A}_N(x)-\mathcal{B}_N(x)\|_p$ is bounded above by
\[
	\sum_{i=1}^N \left\| \left[\int_{Q(x, t_{i+1}, t_i)} 
	\sigma(u(s\,,y))\,\xi(\d s\,\d y)\right]^2 -  
	\left[\int_{Q(x, t_{i+1}, t_i)}\sigma\left(u(r_i(y)\,, y) \right) \,\xi(\d s\,\d y) 
	\right]^2 \right\|_p \leq \sum_{i=1}^N I_i J_i,
\]
where
\begin{align*}
	I_i &:= \left\| \int_{Q(x, t_{i+1}, t_i)} 
		\left[ \sigma(u(s\,,y))-\sigma(u(r_i(y)\,,y)) \right]\xi(\d s\,\d y) 
		\right\|_{2p},\text{ and}\\
	J_i &:= \left\| \int_{Q(x, t_{i+1}, t_i)} 
		\left[ \sigma(u(s\,,y))+\sigma(u(r_i(y)\,,y)) \right]\xi(\d s\, \d y)  \right\|_{2p}.
\end{align*}

We estimate $I_i$ first as follows: Because $|s-r_i(y)|\leq (t_{i+1}-t_i)$
uniformly for all $(s\,,y)\in Q(x\,, t_{i+1}\,, t_i)$, 
\eqref{thm:Holder continuity}, the
Burkholder--Davis--Gundy inequality, and Minkowski's inequality
together yield
\[
	I_i^2 \lesssim \int_{Q(x, t_{i+1}, t_i)} \left\| 
	\sigma(u(s\,,y))-\sigma(u(r_i(y)\,,y))\right\|^2_{2p} \,\d s\,\d y
	\lesssim \int_{Q(x, t_{i+1}, t_i)} (t_{i+1}-t_i) \,\d s\,\d y,
\]
uniformly for all $i=1,\ldots,N$. Thus,
$I_i \lesssim t_{i+1}-t_i$ uniformly for all $i=1,\ldots,N$.
Similarly,
\[
	J_i^2  \leq \int_{Q(x, t_{i+1}, t_i)} \left\| \sigma(u(s\,,y)) + 
	\sigma(u(r_i(y)\,,y))\right\|^2_{2p}\,\d s\,\d y\\
	\lesssim t_{i+1}-t_i,
\]
uniformly for all $i=1,\ldots,N$. These estimates of $I_i$ and $J_i$
can be combined to yield
\[
	\|\mathcal{A}_N(x)-\mathcal{B}_N(x)\|_p\lesssim
	\sum_{i=1}^N  (t_{i+1}-t_i)^{3/2}\lesssim N^{-1/2}\to0
	\qquad\text{as $N \to \infty$.}
\]
This concludes Step 1.\\

\emph{Step 2.} 
Define
\[
	\mathcal{C}_N(x) :=\sum_{i=1}^N \int_{Q(x, t_{i+1}, t_i)} \left[
	\sigma(u(r_i(y)\,,y))\right]^2\,\d s\,\d y.
\]
We plan to establish next that $\mathcal{B}_N(x)-\mathcal{C}_N(x)\to0$ in $L^p(\Omega)$
as $N\to\infty$. According to the general theory of Dalang \cite{Dalang}, 
\begin{equation}\label{moments}
	\sup_{s\in[0,t]}\sup_{x\in\R}\EE\left( |u(s\,,x)|^q\right)<\infty
	\qquad\text{for all $q\ge2$.}
\end{equation}
Therefore, the Burkholder--Davis--Gundy inequality 
ensures that the stochastic processes
$\{\mathcal{B}_N(x)\}_{N\ge1}$ and $\{\mathcal{C}_N(x)\}_{N\ge1}$ are both bounded in 
$L^q(\Omega)$ for all $q\geq 2$. Consequently, it suffices to show that 
$\mathcal{B}_N(x)-\mathcal{C}_N(x)\to 0$ in $L^2(\Omega)$ as $N\to\infty$. With this aim in mind,
write
\[
	\EE\left( \left| \mathcal{B}_N(x)-\mathcal{C}_N(x)\right|^2\right) = B_1 + B_2 - 2B_3,
\]
where
\begin{align*}
	B_1 &:=\EE \left[ \left( \sum_{i=1}^N \left|\int_{Q(x, t_{i+1}, t_i)} 
		\sigma(u(r_i(y)\,,y))\,\xi(\d s\,\d y)\right|^2\right)^2\right],\\
	B_2 &:= \EE\left[ \left( \sum_{i=1}^N 
		\int_{Q(x, t_{i+1}, t_i)} \left[\sigma(u(r_i(y)\,, y))\right]^2
		\,\d s\,\d y \right)^2\right],\quad\text{and}\\
	B_3 &:= \EE \left( \sum_{i=1}^N \left| 
		\int_{Q(x, t_{i+1}, t_i)} \sigma(u(r_i(y)\,,y))\,\xi(\d s\,\d y)
		\right|^2 \times \sum_{i=1}^N  \int_{Q(x, t_{i+1}, t_i)}
		\left[\sigma(u(r_i(y)\,,y))\right]^2\,\d s\,\d y \right).
\end{align*}
In order to simplify the notation, define for all $i=1,\ldots,N$,
\[
	 Q_i := \int_{Q(x, t_{i+1}, t_i)} \sigma(u(r_i(y)\,, y)) \,\xi(\d s\,\d y),\quad
	\widetilde{Q}_i  :=
	\int_{Q(x, t_{i+1}, t_i)} \left[ \sigma(u(r_i(y)\,, y))\right]^2 \,\d s\,\d y.
\]
Whenever $j<i$,
\[
	\EE \left[ Q_i^2  Q_j ^2\right] =\EE \left[ \EE \left( Q_i^2  Q_j ^2 
	\mid \mathcal{G}_{t_i, x}\right)\right]
	= \EE \left[Q_j ^2 \EE \left( Q_i^2  \mid \mathcal{G}_{t_i, x}\right)\right) 
	=\EE \left[ \widetilde{Q}_i Q_j^2\right].
\]
Consequently,
\[
	B_1= \EE\mathop{\sum\sum}\limits_{1\le i,j\le N}  Q_i^2 Q_j^2 
	= \EE \sum_{i=1}^N Q_i^4 + 2 \EE\mathop{\sum\sum}\limits_{1\le j < i\le N} 
	\widetilde{Q}_i Q_j^2.
\]
The same conditioning technique yields
\[
	B_3 = \EE \mathop{\sum\sum}\limits_{1\le j < i\le N} 
	Q_j^2 \widetilde{Q}_i + 
	\EE \mathop{\sum\sum}\limits_{1\le j < i\le N} 
	\widetilde{Q}_j \widetilde{Q}_i + \EE \sum_{i=1}^N Q_i^2 \widetilde{Q}_i.
\]
As a result, it follows that
\[
	\EE \left( |\mathcal{B}_N(x) - \mathcal{C}_N(x)|^2\right)
	=\EE \sum_{i=1}^N Q_i^4 
	+ \EE\sum_{i=1}^N \widetilde{Q}_i^2
	- 2 \EE\sum_{i=1}^NQ_i^2 \widetilde{Q}_i.
\]
Thanks to the uniform boundedness of the moments, and
the Burkholder--Davis--Gundy inequality, 
\[
	\EE  ( |\mathcal{B}_N(x) - \mathcal{C}_N(x)|^2 )
	\lesssim\sum_{i=1}^N (t_{i+1}-t_i)^2
	\lesssim N^{-1}\to0\qquad\text{as $N \to \infty$.}
\]
This concludes Step 2.\\

\emph{Step 3.} We are ready to verify \eqref{eq:time quadratic}. Define
\[
	\mathcal{D}(x) := \int_{Q(x,t)} \left[\sigma(u(s-|x-y|\,,y))\right]^2\,\d s\,\d y
	=\sum_{i=1}^N \int_{Q(x,t_{i+1},t_i)}
	[ \sigma(u(s-|x-y|\,,y))]^2\,\d s\,\d y.
\]
In light of Steps 1 and 2, it remains to prove that $\mathcal{C}_N(x)\to\mathcal{D}(x)$ 
in $L^p(\Omega)$ as $N \to \infty$. One can recycle the argument of Step 1 in order to 
the uniform-in-$N$ bound,
\[
	\|\mathcal{C}_N(x)-\mathcal{D}(x)\|_p \lesssim \sum_{i=1}^N (t_{i+1}-t_i)^{3/2}
	\lesssim N^{-1/2}\to0\qquad\text{as $N \to \infty$}.
\]
This completes Step 3,
and hence the proof of \eqref{eq:time quadratic}. \\

{\it Step 4}. (Sketch) In this step we outline the proof of the remaining assertion
\eqref{eq:space quadratic} of Theorem \ref{thm:quadratic}. The details require
considerably-more space, yet not many more ideas, than those in Steps 1--3.

For every $1\leq i \leq N$ and for each $x_i$, let
\[
	L_i (t) := Q(x_i\,, t)\setminus Q(x_{i+1}\,, t)
	\quad\text{and}\quad
	R_i(t):=Q(x_{i+1}\,, t)\setminus Q(x_i\,, t).
\]
By \eqref{mild},
\begin{equation}
	u(t\,,x_{i+1})-u(t\,,x_i) = \int_{R_i(t)} \sigma(u(s\,,y))\,\xi(\d s\,\d y)
	-\int_{L_i(t)}\sigma(u(s\,,y))\,\xi(\d s\,\d y)\,.
\end{equation}
Define
\begin{align*}
	A_N(t) &:=\sum_{i=1}^N \left[ u(t\,,x_{i+1})-u(t\,,x_i) \right]^2,\\
	B_N(t) &:=\sum_{i=1}^N \left[\int_{R_i(t)}\sigma(u(v_i(y)\,, y))\,\xi(\d s\,\d y)-
		\int_{L_i(t)}\sigma(u(v_i(y)\,, y))\,\xi(\d s\,\d y)\right]^2,\\
	C_N(t) &:= \sum_{i=1}^N \int_{L_i(t)\cup R_i(t)} \left[\sigma(u(v_i(y)\,, y))\right]^2
		\,\d s\,\d y,\text{ and}\\
	D(t) &:=  \int_{(0,t)\times(X_1,X_2)} 
		\left(\left[\sigma(u(s\,,x-t+s))\right]^2 + 
		\left[\sigma (u(s\,,x+t-s))\right]^2\right)\d s\,\d x,
\end{align*}
where $v_i(y):= \max(t+y-x_{i+1}\,, 0)$ if $(s\,,y)\in L_i(t)$, 
and $v_i(y):= \max(t-y+x_i\,,0)$ if $(s\,,y)\in R_i(t)$. It is possible to adapt the arguments
of Steps 1--3 in order to prove that 
$A_N(t)-B_N(t)\to0$ in $L^p(\Omega)$ as $N\to\infty$. 
Next we argue that $B_N(t)-C_N(t)\to0$ in $L^p(\Omega)$ as $N\to\infty$. 
It is easy to see that both $\{B_N(t)\}_{N\ge1}$ and $\{C_N(t)\}_{N\ge1}$ 
are uniformly bounded in $L^p(\Omega)$. Therefore,
it suffices to prove the convergence of $B_N(t)-C_N(t)$ to zero in $L^2(\Omega)$. 

In order to save on typography, define
\begin{align*}
	&R_i := \int_{R_i(t)} \sigma(u(v_i(y)\,,y))\,\xi(\d s\,\d y)\,, \quad 
		\hat {R}_{i}:= \int_{R_i(t)} \left[\sigma(u(v_i(y)\,,y))\right]^2\,\d s\,\d y\,, \\
	&L_i:=\int_{L_i(t)} \sigma(u(v_i(y)\,,y))\,\xi(\d s\,\d y)\,,  \quad 
		\hat{L}_{i}: = \int_{L_i(t)} \left[\sigma(u(v_i(y)\,,y))\right]^2\,\d s\,\d y.
\end{align*}
An expansion of the square yields
\[
	\EE \left( |B_N(t)-C_N(t)|^2\right) 
	= \EE\left( \left|\sum_{i=1}^N (R_i-L_i)^2-\sum_{i=1}^N(\hat{R}_i 
	+ \hat{L}_i) \right|^2\right)
	:= S_1 + S_2 -2S_3,
\]
where
\[
	S_1:=\EE \left(\left|\sum_{i=1}^N (R_i-L_i)^2 \right|^2\right)
	\quad\text{and}\quad
	S_2 := \EE\left(\left| \sum_{i=1}^N (\hat{R}_i + \hat{L}_i) \right|^2\right),
\]
and $S_3$ is the remainder.
We compute $S_1$, $S_2$, and $S_3$ in this order.

Let us introduce $\sigma$-algebras $\mathcal{F}_\pm(x_i)$ as follows: 
Let $\mathcal{F}_+(x_i\,,\varepsilon)$ denote the $\sigma$-algebra generated by 
$\int\phi\,\d\xi$ for all smooth functions $\phi$ that are
supported in $\cup_{x_1 \leq x \leq x_i + \varepsilon} Q(x, t)$.
Similarly, let $\mathcal{F}_-(x_i\,,\varepsilon)$ denote the 
$\sigma$-algebra generated by $\int\phi\,\d\xi$ for all smooth
$\phi$ supported on $\cup_{x_i-\varepsilon \leq x \leq x_N} Q(x\,, t)$.
Then, we define
$\mathcal{F}_\pm(x_i) := \cap _{\varepsilon >0} \mathcal{F}_\pm (x_i\,, \varepsilon)$,

If  $1\le j<i\le N$, 
then we may condition on $\mathcal{F}_+(x_i)$ and/or $\mathcal{F}_-(x_{j+1})$ 
in order to see that 
\[
	\EE \left[(R_i-L_i)^2(R_j-L_j)^2\right] = 
	\EE \left[R_i^2 R_j^2 + L_i^2 R_j^2 + R_i^2L_j^2+ L_i^2 L_j^2 \right].
\]
Similar considerations show that the above also holds when $1\le i<j\le N$. Thus, 
\[
	S_1 = \EE 
	\mathop{\sum\sum}\limits_{1\le i\neq j\le N} 
	\left(R_i^2 R_j^2 + L_i^2 R_j^2 + R_i^2L_j^2+ L_i^2 L_j^2 \right) + \EE\sum_{i=1}^N (R_i-L_i)^4.
\]
If $1\le j < i\le N$, another conditioning {argument} yields
\[
	\EE\left[R_i^2 R_j^2 + L_i^2 R_j^2 + R_i^2L_j^2+ L_i^2 L_j^2 \right] = 
	\EE \left[ R_j^2 \hat{R}_i + L_i^2 R_j^2 + \hat{R}_i \hat{L}_j + L^2_i \hat{L}_j \right].
\]
By comparison, if $1\le i<j\le N$,
\[
	\EE\left[R_i^2 R_j^2 + L_i^2 R_j^2 + R_i^2L_j^2+ L_i^2 L_j^2 \right]
	= \EE \left[ R_i^2 \hat{R}_j + \hat{L}_i \hat{R}_j + R_i^2 L_j^2 + \hat{L}_i L_j^2 \right].
\]
Thus, we can rearrange the sum to see that
\[
	S_1 = 2\EE\mathop{\sum\sum}\limits_{1\le j < i\le N} 
	\left( \hat{R}_i R_j^2 + \hat{L}_j \hat{R}_i + L^2_i R^2_j + L^2_i \hat{L}_j \right) + 
	\EE\sum_{i=1}^N (R_i-L_i)^4\,. 
\]

For $S_2$, there is no need for conditioning arguments, as a direct calculation yields
\[
	S_2 = 2 \EE\mathop{\sum\sum}\limits_{1\le j < i\le N} 
	\left( \hat{R}_i \hat{R}_j + \hat{L}_i \hat{R}_j + \hat{R}_i \hat{L}_j + \hat{L}_i \hat{L}_j\right) 
	+ \EE\sum_{i=1}^N (\hat{R}_i + \hat{L}_i)^2\,.
\]

Finally, another conditioning argument shows that 
\begin{align*}
	S_3&=\EE\mathop{\sum\sum}\limits_{1\le j < i\le N} 
		\left( \hat{R}_i \hat{R}_j + L^2_i \hat{R}_j + \hat{R}_i \hat{L}_j + 
		L^2_i \hat{L}_j + R^2_j \hat{R}_i + \hat{L}_j \hat{R}_i  + R^2_j \hat{L}_i + \hat{L}_i \hat{L}_j\right)\\
	&\hskip3in+\EE\sum_{i=1}^N (R_i-L_i)^2 (\hat{R}_i + \hat{L}_i).
\end{align*}
One can combine the preceding and compute to see, after a few lines, that
\begin{align*}
	S_1+S_2-2S_3&=2  \EE \mathop{\sum\sum}\limits_{1\le j < i\le N} 
		( L^2_i - \hat{L}_i) ( R_j^2-\hat{R}_j) 
		+ \EE \sum_{i=1}^N (R_i-L_i)^4 + \EE \sum_{i=1}^N (\hat{R}_i + \hat{L}_i)^2 \\
	&\hskip3in- 2 \EE \sum_{i=1}^N (R_i-L_i)^2 (\hat{R}_i + \hat{L}_i)\,.
\end{align*}
In the cases that $L_i(t)$ does not intersect with $R_j(t)$, it is easy to see that 
\[
	\EE \left[( L^2_i - \hat{L}_i )  ( R_j^2-\hat{R}_j )\right]=0.
\]
In the cases that $L_i(t)$ intersects with $R_j(t)$, let $Q(i\,,j):= L_i(t) \cap R_j(t)$, 
and let $L_{i,1}$ and $R_{j,1}$ respectively denote the parts of $L_i(t)$ and 
$R_j(t)$ that lie above $Q(i\,,j)$. Also, let $L_{i,2}$ and $R_{j,2}$ 
respectively denote the parts of $L_i(t)$ and $R_j(t)$ that lie below $Q(i\,,j)$. 
We can condition on $\mathcal{G}_{t, x_{j+1}} \vee \mathcal{G}_{t, x_{i+1}}$
to see  that 
\begin{align*}
	&\EE \left[( L^2_i - \hat{L}_i )  ( R_j^2-\hat{R}_j )\right]\\
	&= \EE \Bigg[2\int_{L_{i,2}} \sigma(u(v_i(y)\,,y))\,\xi(\d s\,\d y) 
		\int_{Q(i,j)} \sigma(u(v_i(y)\,,y))\,\xi(\d s\,\d y) \\
		&\hskip1in+ \left|\int_{L_{i,2}} \sigma(u(v_i(y),y))\,\xi(\d s\,\d y)\right|^2 
		+\left| \int_{Q(i,j)} \sigma(u(v_i(y),y))\,\xi(\d s\,\d y)\right|^2 \\
	&\hskip1in - \int_{L_{i,2}} \left[\sigma(u(v_i(y)\,,y))\right]^2\,\d s\,\d y - 
		\int_{Q(i,j)} \left[\sigma(u(v_i(y)\,,y))\right]^2\,\d s\,\d y\Bigg] \times\\
	&\quad 
		\times \left[ \left|\int_{R_j(t)} \sigma(u(v_i(y)\,,y))\,\xi(\d s\,\d y)\right|^2 
		- \int_{R_j(t)} \left[\sigma(u(v_i(y)\,,y))\right]^2\,\d s\,\d y \right].
\end{align*}
A few more rounds of conditioning on $\mathcal{G}_{t, x_j}\vee \mathcal{G}_{t, x_i}$ yield
\begin{align*}
	\EE \left[(L^2_i - \hat{L}_i )  ( R_j^2-\hat{R}_j )\right]
	:= \EE \left[(\mathcal{L}_{i,1} + \mathcal{L}_{i,2})
	(\mathcal{R}_{j,1} + \mathcal{R}_{j,2})\right],
\end{align*}
where
\begin{align*}
	 \mathcal{L}_{i,1} :=&  2\int_{L_{i,2}} \sigma(u(v_i(y)\,,y))\,\xi(\d s\,\d y) 
	 	\int_{Q(i,j)}\sigma (u(v_i(y)\,,y))\,\xi(\d s\,\d y) \\
	& +\left| \int_{Q(i,j)} \sigma(u(v_i(y)\,,y))\,\xi(\d s\,\d y)\right|^2  
		- \int_{Q(i,j)} \left[\sigma(u(v_i(y)\,,y))\right]^2\,\d s\,\d y, \\
	 \mathcal{L}_{i,2}:=&\left| \int_{L_{i,2}} \sigma(u(v_i(y)\,,y))\,\xi(\d s\,\d y)\right|^2 
	 	- \int_{L_{i,2}} \left[\sigma(u(v_i(y)\,,y))\right]^2\,\d s\,\d y, \\
	\mathcal{R}_{j,1}:=& 2\int_{R_{j,2}} \sigma(u(v_j(y)\,,y))\,\xi(\d s\,\d y) 
		\int_{Q(i,j)} \sigma(u(v_j(y)\,,y))\,\xi(\d s\,\d y) \\
	&\quad +\left| \int_{Q(i,j)} \sigma(u(v_j(y)\,,y))\,\xi(\d s\,\d y)\right|^2 
		- \int_{Q(i,j)} \left[ \sigma(u(v_j(y)\,,y))\right]^2\,\d s\,\d y,\text{ and}\\
	 \mathcal{R}_{j,2}:=& \left|\int_{R_{j,2}} \sigma(u(v_j(y)\,,y))\,\xi(\d s\,\d y)\right|^2 
	 	- \int_{R_{j,2}} \left[ \sigma(u(v_j(y)\,,y)) \right]^2\,\d s\,\d y.
\end{align*}
We may condition on $\mathcal{G}_{t, x_{i+1}}\cap \mathcal{G}_{t, x_j}$
in order to see that $\EE [\mathcal{L}_{i,2}\mathcal{R}_{j,2}]=0$. 
Therefore, the Minkowski and Cauchy-Schwarz inequalities together show that
\[
	\EE (\mathcal{L}_{i,1} \mathcal{R}_{j,1} +  \mathcal{L}_{i,1} \mathcal{R}_{j,2} 
	+ \mathcal{L}_{i,2}\mathcal{R}_{j,1})\lesssim (x_{i+1}-x_i)^{3/2} (x_{j+1}-x_j)
	\lesssim N^{-5/2},
\]
and hence
\[
	\mathop{\sum\sum}\limits_{1\le j < i\le N} 
	\EE (\mathcal{L}_{i,1} \mathcal{R}_{j,1} +  
	\mathcal{L}_{i,1} \mathcal{R}_{j,2} + \mathcal{L}_{i,2}\mathcal{R}_{j,1})
	\lesssim N^{-1/2}\to 0\qquad\text{as $N\to\infty$.}
\]

In like manner, the Minkowski and Cauchy-Schwarz inequalities also show that 
\[
	\EE \sum_{i=1}^N (R_i-L_i)^4 + \EE \sum_{i=1}^N (\hat{R}_i + \hat{L}_i)^2 
	- 2 \EE \sum_{i=1}^N (R_i-L_i)^2 (\hat{R}_i + \hat{L}_i)\to0\quad\text{as $N\to\infty$}.
\]
Thus, $\lim_{N \to \infty}\EE(|B_N(t)-C_N(t)|^2)=0$, as was announced.
Finally, it is possible to reuse the arguments of Steps 1--3 in order to
show that $C_N(t)\to D(t)$ in $L^p$  as $N\to\infty$. It was shown earlier in
Step 4 that $C_N(t)-B_N(t)\to 0$ and $B_N(t)-A_N(t)\to 0$. Thus,
$A_N(t)\to D(t)$ in $L^p(\Omega)$ as $N\to\infty$, as was desired.\qed

\section{Sketch of the Proof of Theorem \ref{thm:CLT and LIL}}%
	\label{sec:Proof of thm CLT and LIL}

As in the proof of Theorem \ref{thm:quadratic},
define $r(y):= \max(t-|x-y|\,, 0)$, and set
\[
	\widetilde{Q}(t\,, t+\varepsilon\,, x): = 
	\left\{(s\,,y):\, (s\,,y)\in Q(t\,, t+\varepsilon\,,x)\  
	\text{and} \ |y-x|\leq t\right\},
\]
where $Q(t\,,t+\varepsilon\,,x)$ was defined in \eqref{rQQ}.
Next, define for all $h>0$,
\[
	M_h = M_h(t\,,x):= \int_{\widetilde{Q}(t, t+h, x)} \sigma(u(r(y)\,,y))\,\xi(\d s\,\d y),
\]
and set $M_0:=\lim_{h\downarrow0}M_h=0$.
The elementary properties of the Walsh stochastic integral imply 
that, given $\mathcal{G}_{t,x}$,   the process
$\{M_h\}_{h\ge0}$ is conditionally a mean-zero, continuous $L^2(\Omega)$-martingale with 
quadratic variation $\langle M\rangle_h = h\mathcal{V}$, 
where $\mathcal{V}=\mathcal{V}(t\,,x)$ is the 
$\mathcal{G}_{t,x}$-measurable random variable,
\[
	\mathcal{V} := \int_{x-t}^{x+t} \left[ \sigma(u(t-|x-y|\,, y))\right]^2\d y.
\]
Thus, L\'evy's characterization theorem of Brownian motion implies that,
given $\mathcal{G}_{t,x}$, $M$ is conditionally a Brownian motion with
variance $h\mathcal{V}$ at time $h>0$.
As in the proof of Theorem \ref{thm:quadratic} (see Steps 1 and 2 of that proof), 
one can prove that
\begin{equation*}
	\lim_{h\downarrow0}
	\frac{u(t+h\,, x) - u(t\,,x) - M_h}{\sqrt h}
	= 0\qquad\text{in $L^2(\Omega)$}.
\end{equation*}
We omit the details. Instead, we mention only that
the central limit theorem \eqref{eq:CLT} follows immediately from this and
the scaling properties of the [conditional] Brownian motion $M$.

In order to prove the more delicate law of the iterated logarithm of the theorem,
define 
\[
	R_h = R_h(t\,,x):=u(t+h\,, x)-u(t\,,x) - M_h\qquad\text{for all $h\ge0$.}
\]
It is not hard to use \eqref{mild} together with
the Burkholder--Davis--Gundy inequality and \eqref{thm:Holder continuity},
as well as \eqref{moments}, in order to see that, for every $p\ge2$,
$\|R_h\|_p \lesssim h$ uniformly for all $h\in[0\,,1]$.
Since $R_0=0$, the Kolmogorov continuity theorem implies that,
for every fixed $\eta\in(0\,,1)$,
\[
	\left\|\sup_{s\in[0,h]} |R_s|\right\|_p\lesssim h^\eta\qquad\text{%
	uniformly for all $h\in[0\,,1]$.}
\]
Thus, a standard application of the Borel--Cantelli lemma yields
$R_h=o(h^\eta)$ almost surely as $h\downarrow0$. In particular,\
\[
	\limsup_{h \to 0} \frac{u(t+h\,, x) - u(t\,,x)}{\sqrt{2h \log\log(1/h)}} 
	= \limsup_{h \to 0} \frac{M_h}{\sqrt{2h \log\log(1/h)}}=\mathcal{V}^{1/2}
	\quad\text{a.s.}
\]
thanks to Khintchine's LIL for the [conditional]
Brownian motion $M$. The law of the iterated logarithm of the theorem---see
\eqref{eq:LIL}---follows.\qed

\begin{small}

\end{small}


\begin{thebibliography}{999}

\bibitem{Dalang} Robert C. Dalang (1999).
	Extending the martingale measure stochastic integral with applications to 
	spatially homogeneous s.p.d.e.'s. 
	{\it Electron.\ J. Probab.}\ {\bf 4}{\it (6)} 29 pp. (electronic).
%
\bibitem{DF} Robert C. Dalang and N. E. Frangos (1998).
	The stochastic wave equation in two spatial dimensions. 
	{\it Ann.\ Probab.}\ {\bf 26}{\it (1)}:187--212. 
%
\bibitem{FKM} Mohammud Foondun,  Davar Khoshnevisan, and
	Pejman Mahboubi (2015).
	Analysis of the gradient of the solution to a stochastic heat equation via fractional Brownian motion. 
	{\it Stoch.\ Partial Differ.\ Equ.\ Anal.\ Comput.} {\bf 3}{\it (2)}:133--158.  
%
\bibitem{Hairer} Martin Hairer (2013). 
	Solving the KPZ equation. 
	{\it Ann.\ of Math.}\ (2) {\bf 178}{\it (2)}:559--664.
%
\bibitem{Hairer1} M. Hairer (2014).
	A theory of regularity structures, {\it Invent.\ Math.}\ {\bf 198}{(2)}:269--504.
%
\bibitem{HP} Martin Hairer and \'Etienne Pardoux (2015).
	A Wong--Zakai theorem for stochastic PDEs. 
	{\it J. Math.\ Soc.\ Japan} {\bf 67}{\it (4)}:1551--1604. 
%
\bibitem{KSXZ} Davar Khoshnevisan, Jason Swanson, Yimin Xiao, and Liang Zhang (2013).
	Weak existence of a solution to a differential equation driven by a very rough fBm.
	Unpublished manuscript. 
	Preprint available at \url{https://arxiv.org/abs/1309.3613}.
%
\bibitem{PZ} Szymon Peszat and Jerzy Zabczyk (2000).
	Nonlinear stochastic wave and heat equations,
	{\it Probab.\ Theory Related Fields} {\bf 116}{\it (3)}:421--443. 
%
\bibitem{Walsh} John Walsh (1986). 
	{\it An Introduction to Stochastic Partial Differential Equations.}
	\'Ecole d'\'et\'e de probabilit\'es de Saint-Flour XIV (1984) 265-439. In:
	Lecture Notes in Math.\ {\bf 1180} Springer, Berlin. 
\end{thebibliography}
\end{document}